\documentclass{amsart}
\usepackage{amsfonts,amssymb,amsmath,amsthm}
\usepackage{url}
\usepackage{enumerate}

\newtheorem{thrm}{Theorem}[section]
\newtheorem{lem}[thrm]{Lemma}
\newtheorem{prop}[thrm]{Proposition}
\newtheorem{cor}[thrm]{Corollary}
\theoremstyle{definition}
\newtheorem{definition}[thrm]{Definition}

\newtheorem{remark}[thrm]{Remark}

\numberwithin{equation}{section}

\newcommand{\Ass}{\operatorname{Ass}}

\newcommand{\grade}{\operatorname{grade}}

\newcommand{\Assh}{\operatorname{Assh}}
\newcommand{\Spec}{\operatorname{Spec}}

\newcommand{\ara}{\operatorname{ara}}

\newcommand{\cd}{\operatorname{cd}}

\newcommand{\Ht}{\operatorname{ht}}
\newcommand{\pd}{\operatorname{pd}}

\newcommand{\Ext}{\operatorname{Ext}}
\newcommand{\Supp}{\operatorname{Supp}}

\newcommand{\Tor}{\operatorname{Tor}}

\newcommand{\Hom}{\operatorname{Hom}}
\newcommand{\Att}{\operatorname{Att}}
\newcommand{\Ann}{\operatorname{Ann}}
\newcommand{\Rad}{\operatorname{Rad}}

\newcommand{\depth}{\operatorname{depth}}

\newcommand{\vpl}{\operatornamewithlimits{\varprojlim}}

\newcommand{\vil}{\operatornamewithlimits{\varinjlim}}

\newcommand{\fm}{\frak{m}}
\newcommand{\fp}{\frak{p}}
\newcommand{\fq}{\frak{q}}
\newcommand{\fa}{\frak{a}}
\newcommand{\fb}{\frak{b}}

\newcommand{\fn}{\frak{n}}

\author{M.Y. Sadeghi, \ \ M. Eghbali, and\ \ Kh. Ahmadi-Amoli}

\email{}

\address{Department of Mathematics, Payame Noor University, Tehran,
 19395-3697, Iran.}
\email{my.sadeghi@pnu.ac.ir}

\address{Department of Mathematics, Tafresh University,  39518 79611,
Tafresh-Iran.} \email{m.eghbali@yahoo.com}

\address{ Department of Mathematics, Payame Noor University, Tehran,
 19395-3697, Iran.}
\email{khahmadi@pnu.ac.ir}


\begin{document}

\title[On top local cohomology modules]{On top local cohomology modules, Matlis duality and tensor products}

\begin{abstract}
Let $\fa$ be an ideal of a local ring $(R, \fm)$ with $c =
\cd(\fa,R)$ the cohomological dimension of $\fa$ in $R$. In the case
that $c=\dim R$, we first give a bound for $\depth D(H^c_{\fa}(R))$,
where $c>2$ and $(R,\fm)$ is complete.
 Later,  $H^c_{\fa}(R) \otimes_R H^c_{\fa}(R)$, $D(H^c_{\fa}(R))
\otimes_R D(H^c_{\fa}(R))$ and $H^c_{\fa}(R) \otimes_R
D(H^c_{\fa}(R))$ are examined. In the case $c=\dim R-1$, the set
$\Att_R H^{c}_{\fa}(R)$ is considered.
\end{abstract} \maketitle

\section{Introduction} \label{sect1}

Throughout this note, we assume $(R,\fm)$ is a $d$-dimensional
commutative Noetherian local ring and $\fa$ is an ideal of $R$. Let
$M$ be an $R$-module. For an integer $i \in \mathbb{Z}$ let
$H^i_{\fa}(M)$ denote the $i$-th local cohomology module of $M$ with
respect to $\fa$ as introduced by Grothendieck (cf. \cite{G} and
\cite{Br-Sh}). The first non-vanishing cohomological degree of the
local cohomology modules $H^i_{\fa}(R)$ is well understood. It is
the maximal length of a regular sequence on $R$ consisting of
elements of $\fa$. The last non-vanishing amount of local cohomology
modules, instead, is more mysterious. It is known as the
cohomological dimension of $\fa$ in $R$, denoted by $\cd(\fa,R)$.
The $\cd(\fa,R)$ carries many information about the ideal itself.
For instance, $\fa$ cannot be generated by elements fewer than
$\cd(\fa,R)$. For the applications of $\cd(\fa,R)$ we refer the
reader to the fruitful papers \cite{Hu-L} and \cite{W}.

It is a well-known fact that $ \cd(\fa,R) \leq d$ (cf. \cite[Theorem
6.1.2]{Br-Sh}). In the case that $\fa=\fm$, the equality $\cd(\fm,R)
= d$ happens. The Hartshorne-Lichtenbaum Vanishing Theorem gives a
characterization for vanishing of $H^d_{\fa}(R)$. In \cite[Theorem
3.3]{E-Sch} the second author and Schenzel  have expressed the
isomorphism $H^d_{\fa}(R) \cong H^d_{\fm
\widehat{R}}(\widehat{R}/J)$ for a certain ideal $J$ of
$\widehat{R}$. It indicates the properties of $H^d_{\fa}(R)$ via
$H^d_{\fm \widehat{R}}(\widehat{R})$, which is more well-known.

Let $E :=E_R(R/\fm)$ be the injective hull of the $R$-module
$R/\fm$. By $D_R(-)$ we denote the Matlis dual functor
$\Hom_R(-,E)$. For brevity, we often write $D(-)$ for $D_R(-)$ when
there is no ambiguity about the ring $R$. We denote by $\pd_R (M)$
resp. $K_{M}$ the projective dimension resp. the canonical module of
a finite $R$-module $M$ (if it exists). In Section $2$, we assume
that $H^d_{\fa}(R) \neq 0$, i.e. $\cd(\fa,R)=d$ and give some
information on $D(H^d_{\fa}(R))$ and the tensor products of these
modules. First, for a complete local ring $(R,\fm)$ we show that
$\depth D(H^d_{\fa}(R)) \geq \min \{2, d\}$ (see Proposition
\ref{3.1} ). Then, we give an upper bound for $\depth (D(H^d_{\fa}
(R)))$:
$$\begin{array}{ll} \   \grade (\fb, D(H^d_{\fa} (R)))+\dim R/\fb \geq \depth (D(H^d_{\fa} (R))) \geq 2,
\end{array}$$
where $\fb$ is an ideal of $R$ of dimension at most $2$ and $d \geq
3$. Hence, it provides conditions in which Cohen-Macaulayness of
$D(H^d_{\fa}(R))$ fails ($d \geq 3$), see  Remark \ref{3.5}. Later,
we examine the modules   $H^d_{\fa}(R) \otimes_R H^d_{\fa}(R)$,
$D(H^d_{\fa}(R)) \otimes_R D(H^d_{\fa}(R))$ and $H^d_{\fa}(R)
\otimes_R D(H^d_{\fa}(R))$. In particular, in the case $R$ is
complete and $\pd_{R/Q} (K_{R/Q}) <\infty$ for a certain ideal $Q$
of $R$, we give the necessary and sufficient conditions for
Cohen-Macaulayness of $D(H^d_{\fa}(R)) \otimes_{R/Q}
D(H^d_{\fa}(R))$ (see Theorem \ref{4.2}).

In general, the  case $\cd(\fa,R) = d-1$ is more complicated than
the first one. On the other hand, $H^{d-1}_{\fa}(R)$ is more
mysterious and its properties are not so much known. The aim of the
Section $3$ is to obtain more information about the structure of
$H^{d-1}_{\fa}(R)$ to shed light on the obscure aspects of this
issue. We examine $\Att_R (H^{d-1}_{\fa}(R))$ the set of attached
prime ideals of $H^{d-1}_{\fa}(R)$.  As a technical tool we use the
concept of colocalization introduced by Richardson in \cite{R}. Our
 main result in this section is
Theorem \ref{2.2}.

\section{Matlis duality and tensor products}

Throughout this section, $(R,\fm)$ is a local ring of dimension $d$
and $\fa$ is an ideal of $R$. Moreover, we assume that $H^d_{\fa}
(R) \neq 0$.


 \begin{lem} \label{pre} Suppose that $R$ is complete.

\begin{itemize}
  \item [(1)] Suppose that $M$ and $N$ are either Artinian or finitely generated $R$-modules.  Then we have the following isomorphism
 $$   \Hom_R (M,N) \cong \Hom_R (D(N),D(M)).$$
  \item [(2)] Let $M$ and $N$ be Artinian $R$-modules. Then we have the following isomorphism
  $$  D(N) \otimes_R M \cong D(\Hom_R (M,N)).$$
\end{itemize}
 \end{lem}

\begin{proof} (1) By virtue of \cite[Theorem 10.2.12]{Br-Sh} one has $N \cong D(D(N))$. Then, by \cite[Theorem 2.11]{Ro} we get the claim as follows:
$$   \Hom_R (M,N)  \cong \Hom_R (M,D(D(N)))
 \cong \Hom_R (D(N),D(M)).$$

(2) It should be noted that, by Matlis duality $D(N)$ is a finitely
generated $R$-modules. By $\Hom-\otimes$-adjointness we can prove
the claim as follows:
\begin{eqnarray*}
  D(N) \otimes_R M &\cong &D(N) \otimes_R D(D(M))\\
   &\cong & D(N) \otimes_R \Hom_R (D(M),E(R/\fm) ) \\
   &\cong& \Hom_R (\Hom_R (D(N),D(M)), E(R/\fm))\\
   & \cong & \Hom_R( \Hom_R (M,N), E(R/\fm)).
\end{eqnarray*}
The third isomorphism follows by \cite[10.2.16]{Br-Sh}.
\end{proof}

\ We need the definition of the canonical module of a finitely
generated $R$-module $M.$ To this end let $R$ be the epimorphic
image of a local Gorenstein ring $(S,\mathfrak{n})$ and $n = \dim
S.$

\begin{definition} \label{2.3}  (A)
For a finitely generated $R$-module $M$ we consider
\[
K(M) = \Hom_R(H^t_{\mathfrak{m}}(M),E_R(R/\mathfrak m)),  \; t =
\dim M.
\]
Because $H^t_{\mathfrak{m}}(M)$ is an Artinian $R$-module it turns
out that $K(M)$ is a
 finitely generated $\hat R$-module. \\
 (B)
This is closely related to the notion of the canonical module $K_M$
of a finitely generated $R$-module $M.$ To this end we have to
assume that $(R,\mathfrak m)$ is the epimorphic image of a local
Gorenstein ring $(S,\mathfrak{n})$ with $n = \dim  S.$ Then $K_M =
\Ext_S^{n-t}(M,S)$ is called the canonical module of $M.$ By virtue
of the Cohen's Structure Theorem and  (A)  there is an isomorphism
$K(M) \simeq K_{\hat{M}}.$ By the Local Duality Theorem (see
\cite{G}) there are the following isomorphisms
\[
\Hom_R(H^t_{\mathfrak{m}}(M), E_R(R/\mathfrak{m})) \simeq K(M)
\simeq K_{\hat{M}}.
\]
Here $\hat{R}$ resp. $\hat{M}$ denote the $\fm$-adic completion of
$R$ resp. of $M$.
\end{definition}

In the following, we summarize a few results about the last
non-vanishing local cohomology module $H^d_{\fa}(R)$. To this end
note that a finite module $M$ over a Noetherian ring $R$ satisfies
Serre's condition $S_n$ if $\depth M_{\fp} \geq \min \{n,\dim
M_{\fp} \}$ for all $\fp \in \Supp M$.

\begin{lem} \label{collect} Let $(R,\fm)$ be a complete local ring.
\begin{enumerate}
\item[(1)] $H^d_{\fa}(R)$ is an Artinian $R$-module.
\item[(2)] $H^d_{\fa}(R) \cong H^d_{\fm}(R/Q_{\fa}(R))$, where
$  Q_{\fa}(R)= \bigcap_i \fq_i$ and $\fq_i$s run over all primary
components of the zero ideal of $R$ with $\dim R/\fq = d$ and $\Rad
(\fa + \fq) = \fm$.
\item[(3)] $ D(H^d_{\fa}(R)) \cong K_{R/Q_{\fa}(R)}$ is a finitely generated $R$-module.
\item[(4)] $\Ann_{R} H^d_{\fa}(R)= Q_{\fa}(R)$.
\item[(5)] $\dim_R D(H^d_{\fa}(R))=d$.
\item[(6)] The module $ D(H^d_{\fa}(R))$ satisfies Serre's condition $S_2$.
\item[(7)] If $R/Q_{\fa}(R)$ is a Cohen-Macaulay ring,  then the module $D(H^d_{\fa}(R))$ is Cohen-Macaulay.
\end{enumerate}
\end{lem}

\begin{proof} For the proof of (1), see  \cite [Theorem 7.1.6]{Br-Sh}.
For the proof of (2), (3), (4) and (7) we refer the reader to
\cite[3.3, 4.1(c), 4.2(a) and 4.3]{E-Sch}. Part (5) follows by (2)
and (4). Finally, (6) is a consequence of (3) and \cite [Lemma
1.9]{Sch1}.
\end{proof}

 We recall
 the $i$th
formal local cohomology of $M$ with respect to $\fa$ as
$\mathfrak{F}^i_{\fa}(M):={\vpl}_nH^i_{\fm}(M/\fa^n M)$. See
\cite{Sch2} and \cite{E1} for more information. It should be noted
that $H^i_{\fa}(M,N)={\vil}_n \Ext^{i}_{R}(M/\fa^n M, N)$ is known
as the $i$th generalized local cohomology module for all $R$-modules
$M$ and $N$. In the case $R=M$, the module $H^i_{\fa}(M,N)$ is the
same as the ordinary local cohomology module $H^i_{\fa}(N)$.

\subsection{On the depth}

Let $(R,\fm)$ be a complete local ring. It is a well-known result
that $D(H^d_{\fa} (R))$ is a Cohen-Macaulay module of dimension $d$
for $d \leq 2$ (see for example Lemma \ref{collect}(6)). In this
subsection, we examine the depth of Matlis
 dual of top local cohomology modules.

\begin{prop} \label{3.1} Let $(R,\fm)$ be a complete local ring. Then
$$\begin{array}{ll} \  \depth (D(H^d_{\fa} (R)) \geq \min \{ d, 2\}.
\end{array}$$
\end{prop}

\begin{proof}
We prove the claim by induction on $d$. By the explanation before
this Proposition, for $d \leq 2$, the claim is clear. Assume that $d
\geq 3$.

First note that by virtue of Lemma \ref {collect}, we have
$H^d_{\fa} (R) \cong H^d_{\fm} (R/Q)$, where $Q:=Q_{\fa}(R)$. By the
definition of $Q$ we may deduce that $\fm \notin \Ass R/Q$ and so
$\depth R/Q>0$. Hence, there exists an element $x \in \fm$ which is
nonzero-divisor on $R/Q$. Accordingly, we have the following long
exact sequence
$$\begin{array}{ll} \  0 \rightarrow D(H^d_{\fm}(R/Q)) \stackrel{x}{\rightarrow}
D(H^d_{\fm}(R/Q)) \rightarrow D(H^{d-1}_{\fm}(R/(x,Q))) \rightarrow
\cdots.
\end{array}$$
Since $\dim (R/(x,Q))=d-1$, by induction hypothesis one can see that

$$\depth D(H^{d-1}_{\fm}(R/(x,Q))\geq 2.$$

Hence, there are $u \in \fm \setminus \bigcup_{\fp \in
\Ass(D(H^d_{\fm}(R/Q))/x D(H^d_{\fm}(R/Q)))} \fp$. Thus, it follows
that $\depth (D(H^d_{\fm}(R/Q))) \geq 2$.
\end{proof}

According to \cite[Corollary 1.1.4]{Hel} the $\grade(\fa,
D(H^d_{\fa}(R)))$ plays a significant role to see whether $\fa$ is a
set-theoretic complete intersection, i.e. $\Ht(\fa)=\ara(\fa)$. Note
that $\ara(\fa)$ stands for the arithmetic rank of the given ideal
$\fa$. For set-theoretic complete intersection ideals and related
concepts see for instance \cite{E2}. In this direction, we try to
find an explicit relation between $\depth D(H^d_{\fa}(R))$ and
$\grade(\fa, D(H^d_{\fa}(R)))$. To this end, we denote by
$\grade_{\fa}(M, N)$ the least integer $i$ such that
$H^i_{\fa}(M,N)$ is nonzero for finitely generated $R$-modules $M$
and $N$.

\begin{lem}  \label{chand ta}
Let $\fa$ be an ideal of a local ring $(R,\fm)$. Let $M$ and $N$ be finitely generated $R$-modules.
Then
\begin{enumerate}
  \item  $\grade_{\fa}(M, N) \geq \grade_{\fa}(M, R) - \pd_R
  N,$
provided that $\pd_R N < \infty$.
  \item  In the case $(R,\fm)$ is factor ring of a Gorenstein
local ring $(S,\fn)$, one has $\dim M/\fa M = \dim S -
\grade_{\fa}(M, S)$.
  \item  $ \grade_{\fa}(M, N) \geq \max \{0, \depth N - \dim M/\fa M \}$, provided that $\pd_R N < \infty$.
\end{enumerate}
\end{lem}

\begin{proof}
(1) We argue by induction on $n := \pd_R N$.

Put $n = 0$ so $N$ is a finite free $R$-module of the form $N \cong
\oplus_{t} R$, where $t$ is an integer. Hence, by \cite[Theorem
3.4.10]{Br-Sh} the isomorphisms
$$\begin{array}{ll} \
H^i_{\fa}(M, N) \cong H^i_{\fa}(M,\oplus_{t} R) \cong \oplus_{t}
H^i_{\fa} (M, R)
\end{array}$$
hold. Clearly $\grade_{\fa}(M, N) \geq \grade_{\fa} (M, R)$.

Next, suppose that $n > 0$ and the claim holds for $n-1$. There is a
finite free $R$-module $F$ and exact sequence
$$\begin{array}{ll} \ 0 \rightarrow  L  \rightarrow F \rightarrow N \rightarrow 0
\end{array}$$
which induces $\pd_R L = n- 1$. By induction hypothesis
$\grade_{\fa}(M, L) \geq \grade_{\fa}(M, R)-n + 1$. If $i <
\grade_{\fa}(M, R)-n$, then it follows from the long exact sequence
$$\begin{array}{ll} \ \cdots \rightarrow  H^i_{\fa} (M, F)  \rightarrow H^i_{\fa} (M, N) \rightarrow H^{i+1}_{\fa} (M, L) \rightarrow \cdots
\end{array}$$
that $H^i_{\fa} (M, N)$ is zero. Hence, $\grade_{\fa}(M, N) \geq
\grade_{\fa}(M, R) - \pd_R N$.

 (2) In the light of \cite[Remark 3.6 and Theorem 4.5]{Sch2} we
have
$$\begin{array}{ll} \
\dim M/\fa M &= \sup \{i : \mathfrak{F}^i_{\fa}(M) \neq 0\} \\&=
\sup \{i : H^{\dim S-i}_{\fa S} (M, S) \neq 0 \}
\\& = \dim S - \inf \{i : H^i_{\fa S}(M, S) \neq 0 \}
\\& = \dim S - \grade_{\fa}(M, S).
\end{array}$$

(3) First note that as local cohomology, depth and dimension of a
module are stable under completion, by passing to the completion, we
may assume that $R$ is a complete local ring. Then, by Cohen's
Structure Theorem, $R$ is a factor ring of a regular local ring $(T,
\fn)$. In the light of the Independence Theorem, we may assume that
$\mathfrak{F}^i_{\fa}(M) \cong \mathfrak{F}^i_{\fb}(M)$, where $\fb
:= \fa \cap T$.

Using part (2), we have
$$\begin{array}{ll} \ \dim M/\fa M = \dim T - \grade_{\fa}(M, T).
\end{array}$$
So by part (1) and the Auslander-Buchsbaum Theorem, we conclude the
assertion as follows:
$$\begin{array}{ll} \ \dim M/\fa M &= \dim T - \grade_{\fa}(M, T)
\\& \geq \dim T - \grade_{\fa}(M, N) - \pd_T N
\\& = \depth N - \grade_{\fa}(M, N).
\end{array}$$
\end{proof}

 Assume that $\fb$ is an ideal of $R$ of dimension at most $ 2$ and $d \geq 3$.
  By virtue of  Proposition \ref{3.1} and Lemma \ref{chand ta} we have
$$ \grade (\fb, D(H^d_{\fa} (R)))+\dim R/\fb \geq \depth (D(H^d_{\fa} (R))) \geq
2,$$ Provided that $\pd(D(H^d_{\fa} (R)))<\infty$. This inequality
yields an upper bound for the $\depth (D(H^d_{\fa} (R)))$.

\begin{remark}  \label{3.5} Let $(R,\fm)$ be a  local ring of dimension $d >0$.

(1) Let $x$ be a member of a system of parameters of $R$. Consider
$0:_R <x>:= \bigcup_{n \geq 1} 0:_R x^n=0:_R x^t$, for some integer
$t$. The following exact sequence
$$\begin{array}{ll} \   0 \rightarrow 0:_R <x> \rightarrow R \rightarrow R/0:_R <x> \rightarrow 0.
\end{array}$$
implies that $\dim (0:_R <x>)<d$ and $\dim R/(0:_R <x>) =d$. To this
end, note that the element $x$ is not contained in any minimal prime
ideal of $R$. Hence $x^t \subseteq (0:_R (0:_R x^t)) \nsubseteqq
\fp$, for any minimal primes $\fp$ of $R$.

Now apply $H^d_{\fa} (-)$ to the above short
 exact sequence to get the following exact sequence
$$\begin{array}{ll} \  H^d_{\fa} ( 0:_R <x>) \rightarrow H^d_{\fa} (R) \rightarrow H^d_{\fa} (R/0:_R <x>) \rightarrow 0.
\end{array}$$
Hence, by Grothendieck's Vanishing Theorem, we have the isomorphism
$H^d_{\fa} (R) \cong H^d_{\fa} (R/0:_R <x>)$.

(2) Now, we give conditions that Cohen-Macaulayness of $D(H^d_{\fa}
(R))$ fails. Let $(R,\fm)$ be a complete local ring, $d \geq 3$ and
$\fb$ be a one-dimensional ideal containing a regular element $x$ of
$R$ (part (1) guarantees the existence of such an element). Then, it
implies the following long exact sequence
$$\begin{array}{ll} \  \cdots \rightarrow H^{d-1}_{\fa} (R/xR) \rightarrow H^d_{\fa} (R) \stackrel{x}{\rightarrow} H^d_{\fa} (R) \rightarrow 0.
\end{array}$$
By applying $D(-)$ we get
$$\begin{array}{ll} \  0 \rightarrow D(H^d_{\fa} (R))  \stackrel{x}{\rightarrow}
 D(H^d_{\fa} (R)) \rightarrow  D(H^{d-1}_{\fa} (R/xR)) \rightarrow \cdots.
\end{array}$$
Now, it yields that $x$ is a $D(H^d_{\fa} (R))$-regular element. So
$\grade (\fb, D(H^d_{\fa} (R))) >0$. Suppose that $\fb$ does not
contain a $D(H^d_{\fa} (R))/x D(H^d_{\fa} (R))$-regular element.
Hence, \cite[1.2.3 and 1.2.4]{Bru-Her} imply that $\Ext^1_R (R/\fb
,D(H^d_{\fa} (R) )) \neq 0$. So, it follows that $\grade (\fb,
D(H^d_{\fa} (R))) = 1$. Therefore, in the view of Proposition \ref
{3.1} and Lemma \ref{chand ta} we have
$$\begin{array}{ll} \   2 \leq \depth (D(H^d_{\fa} (R)) \leq \grade (\fb, D(H^d_{\fa} (R)))+\dim R/\fb =2.
\end{array}$$
It shows that with the above assumptions $D(H^d_{\fa} (R))$ is not
Cohen-Macaulay for $d \geq 3$.
\end{remark}

\subsection{Tensor products}

 We assume that $H^d_{\fa}(R) \neq 0$ and $Q:= Q_{\fa}(R)$. In
recent research, there is an interest in endomorphism rings of
certain local cohomology modules $H^i_{\fa}(R).$ This was done in
the case of $i = d$  by the second author and Schenzel (see
\cite{E-Sch}) and some conditions for Cohen-Macaulayness of $\Hom_R
(H^d_{\fa}(R), H^d_{\fa}(R))$ have been investigated. For the sake
of completeness, we state the following result.

\begin{remark}  \label{4.1}

 Let $\fa$ be an ideal of a complete local ring $(R, \fm)$. Then,
$$\begin{array}{ll} \   \depth_R (\Hom_R (H^d_{\fa}(R), H^d_{\fa}(R))) \geq \min \{2, d\}.
\end{array}$$
To this end note that $\Hom_R (H^d_{\fa}(R), H^d_{\fa}(R)) \cong
\Hom_R (D(H^d_{\fa}(R)), D(H^d_{\fa}(R))) $ (see Lemma \ref{pre}).
Hence, the inequality is clear by virtue of  \cite[Exercise
1.4.19]{Bru-Her} and Proposition \ref{3.1}.
\end{remark}

The next part of our work is to investigate the tensor product of
top local cohomology modules and their Matlis duality.

\begin{thrm} \label{4.2} Let $(R,\fm)$ be a complete local ring of dimension $d >0$. Then, the following statements hold:
\begin{enumerate}
\item[(1)] $H^d_{\fa}(R) \otimes_R H^d_{\fa}(R)=0$.

\item[(2)] The ideal $Q$ is contained in $\Ann_R (D(H^d_{\fa}(R)) \otimes_R D(H^d_{\fa}(R)))$ and the Krull dimension of
$D(H^d_{\fa}(R)) \otimes_R D(H^d_{\fa}(R))$ is equal to $d$.

\item[(3)]
 Suppose that $\pd_{R/Q} (K_{R/Q}) <\infty $, then the following statements are equivalent :
\begin{enumerate}
\item[(a)] $R/Q$  is a Cohen-Macaulay ring.
\item[(b)]
$$\begin{array}{ll}
\Tor^{R/Q}_i (D(H^d_{\fa}(R)),D(H^d_{\fa}(R)))=0 \text {\ \ \ for
all \ } i \geq 1
\end{array}$$
and
 $$\begin{array}{ll} D(H^d_{\fa}(R)) \otimes_{R/Q} D(H^d_{\fa}(R))  \text {\ \ \ is a Cohen-Macaulay R/Q-module. \ }
\end{array}$$
\end{enumerate}

\item[(4)] Assume that $R/Q$ satisfies $S_2$ situation. Then, there is an isomorphism
$$\begin{array}{ll}
 H^d_{\fa}(R) \otimes_R D(H^d_{\fa}(R))  \cong E_{R/Q}(R/\fm).
\end{array}$$

\end{enumerate}
\end{thrm}

\begin{proof}
(1) As $H^d_{\fa}(R)$ is an Artinian and $\fm$-torsion $R$-module
(cf. Lemma \ref{collect}), then it follows by  \cite[Proposition
5.8]{K-L-S} that $H^d_{\fa}(R) \otimes_R H^d_{\fa}(R)=0$. To this
end note that from Proposition \ref{3.1}, we have $\depth
(D(H^d_{\fa} (R)) > 0$.

(2) It is clear that
$$\begin{array}{ll}
 \Ann_R (D(H^d_{\fa}(R)) \otimes_R D(H^d_{\fa}(R)) )\supseteq \Ann_R
 (D(H^d_{\fa}(R)))=Q.
\end{array}$$
To this end, note that $\Ann_R
 (D(H^d_{\fa}(R)))=\Ann_R
 (H^d_{\fa}(R))$ (see for instance \cite[Remark 10.2.2]{Br-Sh}) then using Lemma \ref{collect} one has $\Ann_R
 (D(H^d_{\fa}(R)))=Q$.

 On the other hand, the following equalities complete the claim.
$$\begin{array}{ll}
 \dim_R (D(H^d_{\fa}(R)) \otimes_R D(H^d_{\fa}(R))) &= \dim_R (\Supp (D(H^d_{\fa}(R)) \otimes_R D(H^d_{\fa}(R))))
 \\&= \dim_R D(H^d_{\fa}(R))=d
\end{array}$$
  (see \ref{pre}).

(3) It should be noted that for  a finite module $M$ which has
$R$-module
 and $R/Q$-module structures then $\depth_R M=\depth_{R/Q}M$. Moreover,
 by virtue of \ref{collect}, it follows that $D(H^d_{\fa}(R)) \cong K_{R/Q}$.

In order to prove that (a) implies (b), take into account that if
$R/Q$ is a Cohen-Macaulay ring, then $D(H^d_{\fa}(R))$ is a
Cohen-Macaulay module (see \ref{collect}). Thus by the
Auslander-Buchsbaum Theorem, we deduce that $K_{R/Q}$ is a free
$R/Q$-module. Hence,
$$\begin{array}{ll}
 \Tor^{R/Q}_i (D(H^d_{\fa}(R)),D(H^d_{\fa}(R)))=0, \text{\  for all\ } i \geq 1.
\end{array}$$
Now, using \cite[Lemma 2.2]{Hu-W} we have
$$\begin{array}{ll} \   \depth (D(H^d_{\fa}(R)) \otimes_{R/Q} D(H^d_{\fa}(R))) = d.
\end{array}$$
Thus, $D(H^d_{\fa}(R)) \otimes_{R/Q} D(H^d_{\fa}(R))$ is a
Cohen-Macaulay $R/Q$-module.

To show that (b) implies (a), we first use the hypothesis and
\cite[Lemma 2.2]{Hu-W}. Then we have
 \begin{eqnarray}
 \ \ \ \ \ \ \ \ \ \ \ 2 \depth D(H^d_{\fa}(R))= \depth R/Q+ \depth (D(H^d_{\fa}(R)) \otimes_{R/Q} D(H^d_{\fa}(R))). \label{dep}
\end{eqnarray}
 As $\pd_{R/Q} (K_{R/Q})<\infty$, then by the Auslander-Buchsbaum Theorem it follows that
 $$\depth D(H^d_{\fa}(R)) \leq \depth R/Q.$$
Since $D(H^d_{\fa}(R)) \otimes_{R/Q} D(H^d_{\fa}(R))$ is a
Cohen-Macaulay $R/Q$-module, then in the light of the equality
(\ref{dep}), we prove the claim as follows:
$$\begin{array}{ll} \  \depth R/Q & \geq  \depth (D(H^d_{\fa}(R)) \otimes_{R/Q} D(H^d_{\fa}(R)))
\\& = \dim_{R/Q} (D(H^d_{\fa}(R)) \otimes_{R/Q} D(H^d_{\fa}(R)))
\\& =\dim R/Q.
\end{array}$$

(4) As $H^d_{\fa}(R)$ is an Artinian $R$-module so $D(H^d_{\fa}(R))$
is finite. Hence,
$$\begin{array}{ll}
 H^d_{\fa}(R) \otimes_R D(H^d_{\fa}(R)) & \cong D( \Hom_R (H^d_{\fa}(R), H^d_{\fa}(R)))\\& \cong D(R/Q) \\& \cong E_{R/Q}(R/\fm).
\end{array}$$
To this end note that the first isomorphism follows by Lemma
\ref{2.2} and the second one follows by \cite[Theorem 4.2]{E-Sch}.
\end{proof}

Let $M$ and $N$ be nonzero finitely generated modules over a local
ring $(R,\fm)$. $M$ and $N$ satisfy the depth formula provided that
$$\begin{array}{ll} \   \depth M+ \depth N =\depth R + \depth (M \otimes_R N).
\end{array}$$
A necessary condition for the depth formula to hold is that $\depth
M+ \depth N \geq \depth R$ (see for instance \cite[pp. 163]{Hu-W}).
As an application of  \ref{3.1} we deduce the following result:

\begin{cor} \label{4.3} Let $(R,\fm)$ be a $d$-dimensional complete local  ring.
 In the case $d=3 \text{\ or } 4$, then we have
$$\begin{array}{ll} \   \depth (D(H^d_{\fa}(R)) \otimes_R D(H^d_{\fa}(R))) \geq 4-\depth R.
\end{array}$$
\end{cor}

\begin{cor} \label{4.4} Let $(R,\fm)$ be a complete local  ring with $d \geq 4$.
Suppose that $\ \pd_{R/Q} (K_{R/Q}) <\infty $ and $\ \Tor^{R/Q}_i
(D(H^d_{\fa}(R)),D(H^d_{\fa}(R)))=0$ for all $i>0$. Then,
\begin{enumerate}
\item[(1)] $\depth(D(H^d_{\fa}(R)) \otimes_{R/Q} D(H^d_{\fa}(R))) \geq 4- \depth R/Q$,

\item[(2)]
$\pd_{R/Q} (D(H^d_{\fa}(R)) \otimes_{R/Q} D(H^d_{\fa}(R))) \leq 2
\depth R/Q-4$.
\end{enumerate}

\end{cor}

\begin{proof} By Lemma \ref{collect} it follows that $D(H^d_{\fa}(R)) \cong K_{R/Q}$ is a
 finitely generated module. Part (1) is a consequence of \ref{3.1} and \cite[Theorem 1.2]{A}.
 Part (2) follows by (1) and the  Auslander-Buchsbaum Theorem.
\end{proof}

\section{Attached primes of $H^{d-1}_{\fa}(R)$}

Throughout this section, let $(R,\fm)$ be a local ring and $S$ a
multiplicative closed subset of $R$. In his paper \cite{R}, A. S.
Richardson  proposed the definition of colocalization of an
$R$-module $M$ relative to $S$  as the $S^{-1}R$-module
$$S_{-1}M :=
D_{S^{-1}R}(S^{-1}D_R(M)).$$ If $S = R \setminus \fp$ for some $\fp
\in \Spec R$, we write $^{\fp}M$ for $S_{-1}M$. See also \cite{E1}
for more information.

In the light of \cite[Theorem 2.2]{R},  representable modules are
preserved  under colocalization and in the case $M$ is a
representable module, we have
\begin{equation}\label{att}
 \ \Att_{S^{-1}R} S_{-1}M = \{ S^{-1} \fp : \fp \in \Att_R M \text{ and } S \cap \fp = \emptyset
\},
\end{equation}
where the notation $\Att_R M$ is used to denote the set of attached
prime ideals of $M$. This section is based on the use of this
property. We examine the set of attached prime ideals of the last
non-vanishing value of local cohomology.

\begin{prop} \label{A}
Let $\fa$ be an ideal of $R$ and $\fp \in \Spec R$. Let $c$ be an
integer such that  $H^i_{\fa}(R)=0$ for every $i > c$. Assume that
$H^c_{\fa}(R/I)$ is representable for all ideals $I$ of $R$. Then,
we have the following
\begin{itemize}
  \item [(1)] $\Att_{R_{\fp}} (^{\fp} H^c_{\fa}(R))\subseteq \{ \fq R_{\fp}:
\dim R/\fq \geq c,  \fq \subseteq \fp \text{ and } \fq \in \Spec R
\}$.
  \item [(2)] In the case $R$ is complete, we get the following equality
$$ \Att_{R_{\fp}} (^{\fp} H^d_{\fa}(R))=
\{ \fq R_{\fp}: \dim R/\fq  = d \text{ , }  \fq \subseteq \fp , \Rad
(\fa+\fq)= \fm \text{ and } \fq \in \Spec R \}.$$
\end{itemize}
\end{prop}

\begin{proof}

(1)  By virtue of (\ref{att}), we have
$$\begin{array}{ll} \ \Att_{R_{\fp}} (^{\fp} H^c_{\fa}(R))=  \{ \fq R_{\fp}: \fq \in \Att_R (H^c_{\fa}(R)) \text{ and }  \fq \subseteq \fp \}.
\end{array}$$
As the functor $H^c_{\fa}(-)$ is right exact and it preserves direct
limits, $H^c_{\fa}(R/\fq)=H^c_{\fa}(R) \otimes_R R/\fq$. On the
other hand, $H^c_{\fa}(R/\fq)$ is representable, then using
\cite[Lemma 2.11]{A-M} one can see that
$$ \Att_R (H^c_{\fa}(R/\fq))=\Att_R (H^c_{\fa}(R)) \cap
\Supp(R/\fq) \supseteq \fq.$$
 Hence,  $\Att_R
(H^c_{\fa}(R/\fq))\neq \emptyset$ which implies that
$H^c_{\fa}(R/\fq)\neq 0$ and consequently $\dim R/\fq \geq c$.

(2) If $H^d_{\fa}(R)=0$, then we are done. Thus we assume that
$H^d_{\fa}(R)\neq 0$.

$"\subseteq"$: Let $\fq R_{\fp} \in \Att_{R_{\fp}} (^{\fp}
H^d_{\fa}(R))$. As we have seen in part (1), $\fq \subseteq \fp$ and
$H^d_{\fa}(R/\fq) \neq 0$ so $\dim R/\fq = d$. Now, it follows from
the Hartshorne-Lichtenbaum Vanishing Theorem that
$\Rad(\fa+\fq)=\fm$.

$"\supseteq"$: Again using (\ref{att}) we have
$$\begin{array}{ll} \ \Att_{R_{\fp}} (^{\fp} H^d_{\fa}(R))= \{ \fq R_{\fp}: \fq \in \Att_R  H^d_{\fa}(R),  \fq \subseteq \fp \},
\end{array}$$
so it is enough to show that $\fq \in \Att_R  H^d_{\fa}(R)$.

As $\dim R/\fq=d$ and $\Rad (\fa+ \fq)=\fm$, so the Independence
Theorem implies that $H^d_{\fa}(R/\fq)\neq 0$. Hence,
\begin{equation}\label{attint}
\ \emptyset \neq \Att_R (H^d_{\fa}(R/\fq))=  \Att_R (H^d_{\fa}(R))
\cap \Supp_R (R/\fq).
\end{equation}
In the contrary assume that $\fq \notin \Att_R  H^d_{\fa}(R)$. Then
by virtue of (\ref{attint}) there exists a prime ideal $\mathcal{Q}
 \in \Att_R H^d_{\fa}(R)$ such that $\mathcal{Q} \supset \fq$ and so
$\dim R/\mathcal{Q} <d$. On the other $\mathcal{Q} \in \Att_R
H^d_{\fa}(R)$ if and only if $\mathcal{Q} R_{\mathcal{Q}}  \in
\Att_{R_{\mathcal{Q}}} (^{\mathcal{Q}} H^d_{\fa}(R))$. By virtue of
part one $\dim R/\mathcal{Q} = d$ which is a contradiction. Now the
proof is complete.
\end{proof}

    After this preparation we can state the main result of this section
as follows:

\begin{thrm} \label{2.2} Let $(R,\fm)$ be a complete local ring of dimension $d$. Let $\fa$
be an ideal of $R$. Assume that $H^{d-1}_{\fa}(R)$ is representable
and $H^d_{\fa}(R)=0$. Then
\begin{itemize}
  \item [(1)] $\Att_R ( H^{d-1}_{\fa}(R)) \subseteq \{ \fp \in \Spec R : \dim R/\fp  = d-1, \Rad (\fa+\fp)= \fm \} \cup \Assh(R)$.
  \item [(2)] $\{ \fp \in \Spec R : \dim R/\fp  = d-1, \Rad (\fa+\fp)= \fm \} \subseteq \Att_R ( H^{d-1}_{\fa}(R))$.
\end{itemize}
\end{thrm}

\begin{proof}

(1) Let $\fp \in \Att_R  H^{d-1}_{\fa}(R)$, then $\fp R_{\fp} \in
\Att_{R_{\fp}} (^{\fp} H^{d-1}_{\fa}(R))$, (cf. (\ref{att})). It
follows from Proposition \ref{A} that $\dim R/\fp \geq d-1$.

If $\dim R/\fp = d$, then it follows that  $\fp \in \Assh(R)$. In
the case $\dim R/\fp = d-1$,  as  $\fp \in \Att_R  H^{d-1}_{\fa}(R)$
and $H^{d-1}_{\fa}(-)$ is a right exact functor, so one can deduce
that $H^{d-1}_{\fa}(R/\fp) \neq 0$. By the Hartshorne-Lichtenbaum
Vanishing Theorem there exists a prime ideal $\fq \supseteq \fp$ of
$R$ of dimension $d-1$ with $\Rad (\fa+\fq)=\fm$. So, it is clear
that $\fq=\fp$.

(2) Let $ \dim R/\fp  = d-1$ and  $\Rad (\fa+\fp)= \fm$, then
Proposition \ref{A}(2) implies that  $\fp R_{\fp} \in \Att_{R_{\fp}}
(^{\fp} H^{d-1}_{\fa}(R/\fp))$. By (\ref{att})  we deduce that  $\fp
\in \Att_R( H^{d-1}_{\fa}(R/\fp))$. On the other hand the
epimorphism
$$\begin{array}{ll} \  H^{d-1}_{\fa}(R) \rightarrow H^{d-1}_{\fa}(R/\fp) \rightarrow 0
\end{array}$$
implies that  $\fp \in \Att_R( H^{d-1}_{\fa}(R))$.

\end{proof}

It is noteworthy to say that in the situation of Theorem \ref{2.2}
if $\fa$ is an ideal of dimension one, the inclusion at (1) will be
an equality, see  \cite[Theorem 8.2.3]{Hel}.

\end{document}